\documentclass[fleqn]{mat01}
\usepackage{times,mathtimy,amssymb}
\begin{document}

\setcounter{page}{165}
\firstpage{165}

\def\theor{\trivlist \item[\hskip \labelsep{\bf Theorem.}]}

\newtheorem{theo}{\bf Theorem}
\renewcommand\thetheo{\Alph{theo}}

\newtheorem{lem}{Lemma}

\font\xx=msam5 at 10pt
\def\qed{\mbox{\xx{\char'03\!}}}

\font\xx=msam5 at 9pt
\def\ab{\mbox{\xx{\char'03}}}

\font\sa=tibi at 10.4pt

\title{A note on absolute summability factors}

\markboth{H~S~{\"O}zarslan}{A note on absolute summability factors}

\author{H~S~{\"O}ZARSLAN}

\address{Department of Mathematics, Erciyes University,
38039 Kayseri, Turkey\\
\noindent E-mail: seyhan@erciyes.edu.tr}

\volume{113}

\mon{May}

\parts{2}

\Date{MS received 25 March 2002; revised 11 March 2003}

\begin{abstract}
In this paper, by using an almost increasing and
$\delta$-quasi-monotone sequence, a general theorem on
$\varphi-{\mid\!{C},\alpha\!\mid}_k$ summability factors, which generalizes
a result of Bor \cite{3} on ${\varphi-\mid\!{C},1\!\mid}_k$ summability factors,
has been proved under weaker and more general conditions.
\end{abstract}


\maketitle

\section{Introduction}

A sequence $(b_n)$ of positive numbers is said to be
$\delta$-quasi-monotone, if $b_n \rightarrow 0$, $b_n>0$ ultimately and
$\Delta b_n \geq -\delta_n$, where $(\delta_n)$ is a sequence of
positive numbers (see \cite{3}). Let $(\varphi_n)$ be a sequence of complex
numbers and let $\sum{a_n}$ be a given infinite series with partial sums
$(s_n)$. We denote by $\sigma_n^\alpha$ and $t_n^\alpha$ the $n$th Ces{\`
a}ro means of order $\alpha$, with $\alpha>-1$, of the sequences $(s_n)$
and $(na_n)$, respectively, i.e.,
\begin{align}
{\sigma_n^\alpha}&=\frac{1}{A_n^\alpha} \sum_{v=0}^{n}{A_{n-v}^{\alpha-
1}}s_{v},\\[.4pc]
{t_n^\alpha}&=\frac{1}{A_n^\alpha} \sum_{v=1}^{n}{A_{n-v}^{\alpha-
1}}va_{v},
\end{align}
where
\begin{equation}
{A_n^\alpha}=O(n^\alpha), \quad \alpha>-1, \quad {A_0^\alpha=1}
\quad \hbox{and} \quad {A_{-n}^\alpha=0} \quad \hbox{for} \quad  n>0.
\end{equation}
The series $\sum{a_n}$ is said to be summable
${\mid\!{C},\alpha\!\mid}_k$, $k \geq 1$ and $\alpha >-1$, if
(see \cite{6})
\begin{equation}
\sum_{n=1}^{\infty} n^{k-1} \mid\!{\sigma}_n^{\alpha}
-{\sigma}_{n-1}^{\alpha}\!\mid^k =
\sum_{n=1}^{\infty} \frac{1}{n} \mid\!t_n^{\alpha}\!\mid^k < \infty,
\end{equation}
and it is said to be summable
${\mid\!{C},\alpha;\beta\!\mid}_k$, $k \geq 1$, $\alpha >-1$ and $\beta \geq 0$,
if (see \cite{7})
\begin{equation}
\sum_{n=1}^{\infty} n^{\beta k+k-1} \mid\!{\sigma}_n^{\alpha}
-{\sigma}_{n-1}^{\alpha}\!\mid^k =
\sum_{n=1}^{\infty} n^{\beta k-1} \mid\!t_n^{\alpha}\!\mid^k < \infty.
\end{equation}
The series $\sum{a_n}$ is said to be summable
$\varphi-{\mid\!{C},\alpha\!\mid}_k$, $k \geq 1$ and $\alpha>-1$,
if (see \cite{2})
\begin{equation}
\sum_{n=1}^{\infty} n^{-k} \mid\!{\varphi_n} t_n^\alpha\!\mid^k < \infty.
\end{equation}
In the special case when ${\varphi_n}=n^{1-\frac{1}{k}}$ (resp.
${\varphi_n}=n^{\beta+1-\frac{1}{k}})$ $\varphi-{\mid\!{C},\alpha\!\mid}_k$
summability is the same as ${\mid\!{C},\alpha\!\mid}_k$ (resp.
${\mid\!{C},\alpha;\beta\!\mid}_k)$ summability.
Bor \cite{4} has proved the following theorem for $\varphi-{\mid\!{C},1\!\mid}_k$
summability factors of infinite series.

\begin{theo}[\!]
Let $(X_n)$ be a positive non-decreasing sequence and
let $(\lambda_n)$ be a sequence such that
\begin{equation}
\mid\!\lambda_n\!\mid X_n= O(1) \quad \hbox{as} \quad {n\rightarrow\infty},
\end{equation}
\begin{equation}
\sum_{v=1}^{n} v X_v \mid\!{\Delta}^2 \lambda_v\!\mid= O(1) \quad \hbox{as} \quad {n\rightarrow\infty}.
\end{equation}
\end{theo}
If there exists an $\epsilon>0$ such that the sequence
$(n^{\epsilon-k} \mid\!\varphi_n\!\mid^k)$ is non-increasing and
\begin{equation}
\sum_{v=1}^{n} v^{-k} \mid\!{\varphi_v} t_v\!\mid^k =O(X_n)
\quad \hbox{as} \quad {n\rightarrow\infty},
\end{equation}
then the series $\sum{a_n}{\lambda_n}$ is summable
${\varphi -\!\mid\!{C},1\!\mid}_k$, $k \geq 1$.

\section{The main result}

The aim of this paper is to extend Theorem A, by using an almost
increasing and $\delta$-quasi monotone sequence, under weaker and more
general conditions for ${\varphi-\!\mid\!{C},\alpha\!\mid}_k$ summability. For
this we need the concept of almost increasing sequence. A positive
sequence $(b_n)$ is said to be {\it almost increasing} if there exists a
positive increasing sequence ${c_n}$ and two positive constants $A$ and
$B$ such that $ A({c_n}) \leq {b_n} \leq B({c_n})$ (see \cite{1}). Obviously
every increasing sequence is an almost increasing sequence, but the
converse need not be true, as can be seen from the example $b_n=n
{\rm e}^{(-1)^n}$. So we are weakening the hypotheses of the theorem in
replacing the increasing sequence by an almost increasing sequence.

Now, we shall prove the following:

\begin{theor}
{\it Let $(X_n)$ be an almost increasing sequence and the sequence
$(\lambda_n)$ such that condition $(7)$ of Theorem A is satisfied. Suppose
that there exists a sequence of numbers $(A_n)$ such that it is
$\delta$-quasi monotone with $\sum n{A_n}{X_n}$ convergent and $\mid\!{\Delta {\lambda_n}}\!\mid \leq \mid\!A_n\!\mid$ for all n. If there exists
an $\epsilon>0$ such that the sequence $(n^{\epsilon-k}
\mid\!\varphi_n\!\mid^k)$ is non-increasing and if the sequence
$(w_n^\alpha)$, defined by (see {\rm \cite{8}})
\begin{equation}
w_n^\alpha= \left\{ \begin{array}{l@{\qquad}l} \mid\!t_n^\alpha\!\mid, & \mbox{$\alpha=1 $} \\[.5pc]
 \max_{1 \leq v \leq n} \mid\!t_v^\alpha\!\mid, & \mbox{$0 <
\alpha < 1$} \end{array} \right.
\end{equation}
satisfies the condition
\begin{equation}
\sum_{n=1}^{m} n^{-k} (w_n^\alpha \mid\!\varphi_n\!\mid)^k =O(X_m)
\quad {as} \quad {m\rightarrow\infty},
\end{equation}
then the series $\sum{a_n}{\lambda_n}$ is summable
${\varphi-\!\mid\!{C},\alpha\!\mid}_k$, $k \geq 1$,~~$0 < \alpha \leq 1$
and $k {\alpha}+ {\epsilon} >1$.}
\end{theor}
We need the following lemma for the proof of our theorem.

\begin{lem}
{\rm {\rm\cite{5}}}. If $0 < \alpha \leq 1$  and $1 \leq v \leq n$, then
\begin{equation}
\left|\sum_{p=0}^{v}{A_{n-p}^{\alpha-1}}a_{p}\right|
\leq \max_{1 \leq m \leq v} \left| \sum_{p=0}^{m}
{A_{m-p}^{\alpha-1}}a_{p}\right|.
\end{equation}
\end{lem}

\section{Proof of the theorem}

Let $(T_n^\alpha)$, with $0 < \alpha \leq 1$,
be the $n$th $(C,\alpha)$ mean of the sequence $(n{a_n}{\lambda_n})$.
Then, by (2), we have
\begin{equation}
T_n^\alpha = \frac{1}{A_n^\alpha} \sum_{v=1}^{n}{A_{n-v}^{\alpha-1}}
v{a_v}{\lambda_v}.
\end{equation}
Using Abel's transformation, we get
\begin{equation*}
T_n^\alpha = \frac{1}{A_n^\alpha} \sum_{v=1}^{n-1}{\Delta \lambda_v}
\sum_{p=1}^{v}{A_{n-p}^{\alpha-1}}p{a_p}+
\frac{\lambda_n}{A_n^\alpha} \sum_{v=1}^{n}
{A_{n-v}^{\alpha-1}}v{a_v},
\end{equation*}
so that making use of Lemma 1, we have
\begin{align*}
\mid\!T_n^\alpha\!\mid & \leq \frac{1}{A_n^\alpha} \sum_{v=1}^{n-1}
\mid\!{\Delta \lambda_v}\!\mid
\Bigg\vert\!\sum_{p=1}^{v}{A_{n-p}^{\alpha-1}}p{a_p}\Bigg\vert+
\frac{\mid\!\lambda_n\!\mid}{A_n^\alpha}\  \Bigg\vert\sum_{v=1}^{n}
{A_{n-v}^{\alpha-1}}v{a_v}\Bigg\vert\\
&\leq \frac{1}{A_n^\alpha} \sum_{v=1}^{n-1}
{A_v^\alpha}{w_v^\alpha} \mid\!{\Delta \lambda_v}\!\mid
+ \mid\!\lambda_n\!\mid {w_n^\alpha}\\
& = {T_{n,1}^\alpha} + {T_{n,2}^\alpha},
\quad \hbox{say}.
\end{align*}
Since
\begin{equation*}
\mid\!{T_{n,1}^\alpha}+{T_{n,2}^\alpha}\!\mid^k
\leq {2^k} (\mid\!{T_{n,1}^\alpha}\!\mid^k+\mid\!{T_{n,2}^\alpha}\!\mid^k),
\end{equation*}
to complete the proof of the theorem, it is sufficient to show that
\begin{equation*}
\sum_{n=1}^{\infty} n^{-k} \mid\!{\varphi_n}{{T_{n,r}^\alpha}}\!\mid^k
<\infty \quad \hbox{for} \quad r=1,2,
\quad \hbox{by} \quad (6).
\end{equation*}
Now, when $k>1$, applying H{\"o}lder's inequality with indices $k$ and ${k'}$,
where $\frac{1}{k}+\frac{1}{k'}=1$, we get
\begin{align*}
\sum_{n=2}^{m+1} n^{-k} \mid\!{\varphi_n}{{T_{n,1}^\alpha}}\!\mid^ k
&\leq \sum_{n=2}^{m+1} n^{-k}{({A_n^\alpha})^{-k}}\mid\!{\varphi_n}\!\mid^k
\Bigg\{\sum_{v=1}^{n-1}{A_v^\alpha}{w_v^\alpha}\mid\!{\Delta \lambda_v}\!\mid\Bigg\}^k\\[.3pc]
&= O(1) \sum_{n=2}^{m+1} n^{-k}n^{-\alpha k}\mid\!{\varphi_n}\!\mid^k
\Bigg\{\sum_{v=1}^{n-1} v^{\alpha k}{({w_v^\alpha})^{k}} \mid\!{A_v}\!\mid\Bigg\}\\[.3pc]
&\quad\times \Bigg\{ \sum_{v=1}^{n-1}\mid\!{A_v}\!\mid \Bigg\}^{k-1}\\[.3pc]
&= O(1)\sum_{v=1}^{m} v^{\alpha k}{({w_v^\alpha})^{k}} \mid\!{A_v}\!\mid
\sum_{n=v+1}^{m+1} \frac{n^{-k}\mid\!{\varphi_n}\!\mid^k}{n^{\alpha k}}\\[.3pc]
&=O(1)\sum_{v=1}^{m} v^{\alpha k}{({w_v^\alpha})^{k}} \mid\!{A_v}\!\mid
\sum_{n=v+1}^{m+1} \frac{n^{\epsilon-k}\mid\!{\varphi_n}\!\mid^k}
{n^{\alpha k+\epsilon}}\\[.3pc]
&= O(1)\sum_{v=1}^{m} v^{\alpha k}{({w_v^\alpha})^{k}} \mid\!{A_v}\!\mid
v^{\epsilon-k}\mid\!{\varphi_v}\!\mid^k
\sum_{n=v+1}^{m+1} \frac{1}{n^{\alpha k+\epsilon}}
\end{align*}

\begin{align*}
&= O(1)\sum_{v=1}^{m} v^{\alpha k} {({w_v^\alpha})^{k}} \mid\!{A_v}\!\mid
v^{\epsilon-k}\mid\!{\varphi_v}\!\mid^k
\int_{v}^{\infty} \frac{{\rm d}x} {x^{\alpha k +\epsilon}}\\
&= O(1)\sum_{v=1}^{m} v \mid\!{A_v}\!\mid v^{-k}{({w_v^\alpha \mid\!\varphi_v\!\mid})
^{k}}\\
&= O(1)\sum_{v=1}^{m-1} \Delta(v\mid\!{A_v}\!\mid)
\sum_{r=1}^{v} r^{-k} {({w_r^\alpha} \mid\!\varphi_r\!\mid)^{k}}\\
&\quad + O(1)m \mid\!{A_m}\!\mid \sum_{v=1}^{m} v^{-k}
{({w_v^\alpha}\mid\!\varphi_v\!\mid)^{k}}\\
&= O(1)\sum_{v=1}^{m-1} \mid\!\Delta(v\mid\!{A_v}\!\mid)\!\mid {X_v}
+ O(1)m \mid\!{A_m}\!\mid {X_m}\\
&= O(1)\sum_{v=1}^{m-1} v\mid\!{A_v}\!\mid {X_v}
+ O(1) \sum_{v=1}^{m-1} (v+1) \mid\!{A_{v+1}}\!\mid {X_{v}}\\
&\quad + O(1)m \mid\!{A_m}\!\mid {X_m}
= O(1) \quad \hbox{as} \quad {m\rightarrow\infty},
\end{align*}
by virtue of the hypotheses of the Theorem.

Again, since $\mid\!\lambda_n\!\mid=O(1/{X_n})=O(1)$, by (7), we have
\begin{align*}
\sum_{n=1}^{m} n^{-k} \mid\!{\varphi_n}{T_{n,2}^\alpha}\!\mid^ k
&= \sum_{n=1}^{m} \mid\!\lambda_n\!\mid^{k-1} \mid\!\lambda_n\!\mid n^{-k}
{({w_n^\alpha}\mid\!{\varphi_n}\!\mid)^{k}}\\
&= O(1)\sum_{n=1}^{m}  \mid\!\lambda_n\!\mid n^{-k}{({w_n^\alpha}
\mid\!{\varphi_n}\!\mid)^{k}}\\
&= O(1)\sum_{n=1}^{m-1} \mid\!\Delta \mid\!\lambda_n\!\mid\!\mid
\sum_{v=1}^{n} v^{-k}{({w_v^\alpha}\mid\!{\varphi_v}\!\mid)^{k}}\\
&\quad+ O(1) \mid\!\lambda_m\!\mid
\sum_{n=1}^{m} n^{-k}{({w_n^\alpha}\mid\!{\varphi_n}\!\mid)^{k}}\\
&= O(1)\sum_{n=1}^{m-1} \mid\!\Delta \lambda_n\!\mid {X_n}
+ O(1) \mid\!\lambda_m\!\mid {X_m}\\
&= O(1)\sum_{n=1}^{m-1} \mid\!A_n\!\mid {X_n}
+ O(1) \mid\!\lambda_m\!\mid {X_m}\\
&= O(1) \quad {\rm as} \quad {m\rightarrow\infty},
\end{align*}
by virtue of the hypotheses of the Theorem.

Therefore, we get
\begin{equation*}
\sum_{n=1}^{m} n^{-k} \mid\!{\varphi_n} {T_{n,r}^\alpha}\!\mid^k
= O(1) \quad \hbox{as} \quad {m\rightarrow\infty}, \quad \hbox{for} \quad r=1,2.
\end{equation*}
This completes the proof of the Theorem.

If we take $\epsilon=1$ and $\varphi_n=n^{1- \frac{1}{k}}$ (resp.
$\epsilon=1$ and $\varphi_n=n^{\beta+1- \frac{1}{k}}$), then we get a
new result related to ${\mid\!{C},\alpha\!\mid}_k$ (resp.
${\mid\!{C},\alpha;\beta\!\mid}_k$) summability factors.

\end{document}